\documentclass[12pt,a4paper,twoside,final,notitlepage, leqno]{article}
\usepackage[english]{babel}
\usepackage[T1]{fontenc}
\usepackage{epsfig, graphicx, amssymb}
\usepackage{amsmath,amsthm,epsfig,amsfonts}
\usepackage{float}
\usepackage{color}
\def\br {\break}

\newcommand{\monitem}{ \smallskip \noindent $\bullet$ \quad  }
\newcommand{\moneq}{\vspace*{-7pt} \begin{equation} \displaystyle }
\newcommand{\moneqstar}{\vspace*{-6pt} \begin{equation*} \displaystyle }
\newcommand{\monendstar}{\vspace*{-6pt} \end{equation*}   }
\newcommand{\monend}{\vspace*{-7pt} \end{equation}   }
\newcommand{\moneqarraystar}{ \begin{eqnarray*} \displaystyle }
\newcommand{\monendarraystar}{ \end{eqnarray*}   }

\newcommand{\dd}{{\rm d}}
\newcommand{\R}{\mathbb{R}}

\definecolor{vertfonce}{rgb}{0.0, 0.5, 0.0}

\def\section*#1{}

\usepackage{fancyhdr}
\renewcommand{\headrulewidth}{0pt}

\begin{document}

\fancypagestyle{plain}{ \fancyfoot{} \renewcommand{\footrulewidth}{0pt}}
\fancypagestyle{plain}{ \fancyhead{} \renewcommand{\headrulewidth}{0pt}}

~

  \vskip 2.1 cm

\centerline {\bf \LARGE Riemannian formulation of Pontrygin's principle}

\bigskip
\centerline {\bf \LARGE for robotic manipulators}

\bigskip

 \bigskip  \bigskip \bigskip

\centerline { \large    Fran\c{c}ois Dubois$^{ab}$,  Hedy C\'esar Ramírez-de-\'Avila$^{c}$,}

\centerline { \large   Juan Antonio Rojas-Quintero$^{d}$}

\smallskip  \bigskip

\centerline { \it  \small
  $^a$   Conservatoire National des Arts et M\'etiers,}

\centerline { \it  \small Laboratoire de M\'ecanique des Structures et des Systèmes Coupl\'es,   Paris}

\centerline { \it  \small
$^b$   Universit\'e Paris-Saclay, Laboratoire de   Math\'ematiques d'Orsay}

\centerline { \it  \small 
$^b$  Tecnol\'ogico Nacional de M\'exico/I. T. Tijuana, Tijuana 22414, BC, Mexico}

\centerline { \it  \small
$^d$    CONACYT -- Tecnol\'ogico Nacional de M\'exico/I. T. Ensenada, Ensenada 22780, BC, Mexico}


\bigskip  \bigskip

\centerline {18 April 2023 
  {\footnote {\rm  \small $\,$ This contribution will be published
 in the  Proceedings of the Sixteenth International Conference Zaragoza-Pau on Mathematics and its Applications, 
Jaca, 7 - 9 September 2022.}}}



\bigskip  \bigskip
\noindent {\bf \large Abstract}

\noindent
In this work, we consider a mechanical system whose  mass tensor implements 
a scalar product  in a Riemannian manifold.
This system is controlled with the help of forces and torques.
A cost functional is minimized to achieve an optimal trajectory. 
In this contribution, this  cost function  is supposed to be
an arbitrary regular function invariant under a change of coordinates. 
Optimal control evolution based on Pontryagin's principle 
induces a covariant second-order ordinary differential equation
for an adjoint variable featuring the Riemann curvature tensor. 
This second order time evolution is derived in this contribution.

\bigskip   \bigskip   \noindent {\bf \large    Introduction}    

\smallskip \noindent
This work is motivated by the controlled dynamics of articulated systems.
The  Euler-Lagrange equations are classically derived from the knowledge
of kinetic and potential energies. Moreover, the control 
of the system
can be modelled by the addition of external forces and torques.
The search of an optimal dynamics depends on a  
given cost function.
Then Pontryagin's approach \cite{PBGM62} allows the emergence of a control 
law from the minimization of the cost function.
After a remark of Brillouin \cite{Br38}, developed by Lazrak and Vall\'ee \cite {LV95} and 
Rojas-Quintero {\it et al.} 
\cite {Juan-these-2013, RVGSA13, RQ-VC-S-21}: a Riemannian structure is present in such
a system.
With a quadratic cost function, a remarkable result has been obtained in
\cite {Juan-these-2013, RVGSA13, RQ-VC-S-21}: the Lagrange multiplier associated to 
Pontryagin's approach can be interpreted as the forces and torques submitted by the dynamical system.
This property is revisited in this contribution where the  cost function  is not required to be a quadratic function anymore,
but can be taken to be a general nonlinear function instead.

\smallskip \noindent 
In the first section, we clarify  
the previous  choice of a natural Riemann 
metric for robotics. Then in Section~2, we recall very classical results concerning differential operators
on a regular Riemannian manifold. In the next section, the art of derivation suggested by Pontryagin is emphasized.
In section~5, the essential of the work done by one of us \cite {Juan-these-2013}
and published in \cite{DFRV15,RDR22,VRFG13} is briefly presented.
A generalized approach is developped in Section 6: the cost function is no more quadratic as it was in  our  previous works.
Comparing the results for quadratic and general cost functions is emphasized in the conclusion. 

\bigskip   \bigskip   \noindent {\bf \large    1) \quad  Riemannian metric for robotics}    

\smallskip \noindent 
We consider a dynamical system parameterized by a finite number 
of functions of time $ \, q^j(t) $.
The manifold of   states  is denoted by $\, Q $: $ \,  q  \equiv  \{ q^j \}  $. 
In the case of an  articulated system,
the mass metric  $\, M(q) \,$ depends on the general coordinates $\, q^j $. 
This  mass tensor is  symmetric  and positive definite    for each state. 
Then the   kinetic energy
\moneq \label{energie-cinetique} 
K(q, {\dot q})  \equiv  {1\over2} \, \sum_{k \, \ell}  M_{k  \ell} (q) 
\, {\dot q}^k \, {\dot q}^\ell  \,
\monend
is a  positive definite quadratic form 
of the time derivatives $ \, {\dot q}^j $.
The  coefficients $\, M_{k \ell} (q) \, $ are ideal candidates to define 
a Riemannian metric structure  on the configuration space.

\smallskip \noindent
This property has been remarked many years ago by Brillouin \cite{Br38}.
It is also mentioned in the book of Spong and Vidyasagar \cite{SV91}.
In their contribution \cite{LV95}, Lazrak and Vall\'ee emphasize  the tensorial nature of this relation.
From the positivity of the kinetic energy, the mass matrix naturally defines a Riemannian metric.
This fundamental remark is the starting point of our contribution, incorporating Riemannian geometry
in the field of  poly-articulated systems, {\it id est} robotics.

\fancyhead[EC]{\sc{F. Dubois, H. C. Ramírez-de-\'Avila, J. A. Rojas-Quintero}}
\fancyhead[OC]{\sc{Riemannian formulation of Pontrygin's principle}}
\fancyfoot[C]{\oldstylenums{\thepage}}

\bigskip   \bigskip   \noindent {\bf \large    2) \quad  Classical Riemannian geometry}    

\smallskip \noindent
We  follow essentially the presentation of tensorial calculus presented in Lichnerowicz \cite{Li46}.
We use Einstein notation for implicit summation for repeted indices. 
We recall very briefly the main notions.

\smallskip \noindent 
Inverse of the metric mass tensor $ \, M^{-1} $:  $ \, M^{j \ell} $.
We have the contraction   $ \, M_{i j} \, M^{j \ell}  =  \delta^\ell_i \, $
with~$\, \delta^\ell_i \, $ the Kronecker symbol.

\smallskip \noindent 
Covariant space differentiation  
along the manifold
$ \, \partial_j \equiv {{\partial}\over{\partial q^j}} $. The associated 
contravariant basis of the tangent space $\, e_j \, $ is defined by  $\, e_j \equiv \partial_j $. 
The covariant basis $\, e^j \, $  of the tangent space is defined by the relations   
$ \, \, < e^j \,,\,  e_k > = \delta^j_k $,  
where $ \, < . \,,\, .  > \, $ is the duality product between a vector space and its dual.
A contravariant vector field $ \, \varphi = \varphi^k \, e_k \, $ admits also 
covariant components $\, \varphi_j $.  We have the relations
$ \,  \varphi_j = M_{j k} \,  \varphi^k \, $ and conversely $ \,  \varphi^k  =  M^{k j} \,  \varphi_j $ 
between the  contravariant components~$\, \varphi^k \,$ and the covariant components. 

\smallskip \noindent 
Differentiation of a contravariant basis vector 
 $ \, \dd e_j = \Gamma^\ell_{j k} \, \dd q^k \, e_\ell $. It introduces the 
connection $  \, \Gamma^j_{i k} =  {1\over2} \, M^{j \ell} \, \big( \partial_i M_{\ell k} + \partial_k M_{\ell i} - \partial_\ell M_{i  k} \big) $.
These Riemann-Christofell coefficients   $\, \Gamma^j_{k i} \, $  satisfy a symmetry property: 
$  \, \Gamma^j_{k i} = \Gamma^j_{i k} $.
Then the  differentiation of the covariant  basis vector satisfies the relation
$ \, \dd e^j = - \Gamma^j_{k \ell} \, \dd q^k \, e^\ell $.

\smallskip \noindent    
Differentiation of a scalar field $ \, V  \, $: we have $ \, \dd V  = \partial_j V  \, e^j $. Then the 
gradient of the  scalar field~$\, V \ $ satisfies  $ \, \nabla V =   \partial_\ell V \,e^\ell $;
it is a covector field and we have  $ \, {\rm d} V =  \partial_\ell V \, {\rm d}q^\ell$\br
$ = < \nabla V \,,\, {\rm d}q^j \, e_j >  $. 
The covariant  derivative  
of a vector field $ \, \varphi \equiv \varphi^j \, e_j \, $
can be evaluated according to the relation 
   $ \, {\rm d} \varphi  = 
\big( \partial_\ell \varphi^j +  \Gamma^j_{\ell k} \,  \varphi^k  \big) \, {\rm d}q^\ell \, e_j  $.
Analogously, the  covariant  derivative 
 of a    covector field 
$  \xi \equiv \xi_\ell \, e^\ell \, $ satisfies the condition 
$ \,  {\rm d} \xi = \big( \partial_k \xi_\ell -  \Gamma^j_{k \ell} \,  \xi_j  \big) \, {\rm d}q^k \, e^\ell $.
Then the gradient of a  covector field 
satisfies the conditions 
$ \, \nabla \xi = \big( \partial_k \xi_\ell -  \Gamma^j_{k \ell} \,  \xi_j  \big) \, e^k \,  e^\ell $.
It is a two times covariant tensor and  we have 
$  \,   \dd  \xi = \,  < \nabla \xi \,,\,  {\rm d}q^j \, e_j > $. 
Similarly, the second order gradient $ \, \nabla^2 V \, $  
of a scalar field $ \, V \,$ is defined by the relation $ \, \nabla^2 V = \nabla ( \nabla V) $,
{\it id est}  $ \, \nabla^2 V =
\big( \partial_k \partial_\ell V -  \Gamma^j_{k \ell} \,  \partial_j V \big) \, e^k \,  e^\ell $.
It is also a  two times covariant tensor. 

\smallskip \noindent 
Ricci identities  for the differentiation  
of the metric: 
$ \,   \partial_j M_{k \ell}  = \Gamma^p_{j k}  \,  M_{\ell p} +  \Gamma^p_{j \ell}  \,  M_{k  p} $.
We have also 
$ \,   \partial_j M^{k \ell} = - \Gamma^k_{j p} \, M^{p \ell}  - \Gamma^\ell_{j p} \, M^{p k} $. 

\smallskip \noindent
The components $ \,  R^j_{i k \ell} \, $  of the Riemann tensor are defined by the relations  
\moneq \label{riemann-tensor} 
R^j_{i k \ell} \equiv  \partial_\ell \Gamma^j_{i k} - \partial_k \Gamma^j_{i \ell}
+ \Gamma^p_{i k} \, \Gamma^j_{p \ell} - \Gamma^p_{i \ell} \, \Gamma^j_{p k} \, . 
\monend 
We observe the anti-symmetry of the Riemann tensor: 
$ \, R^j_{i k \ell} =  - R^j_{i \ell k} $.
For a given vector field $ \, \varphi \, $ and covector field $ \,\xi $, we introduce the
covector field $ \, R_\varphi . \xi \, $ defined by 
\moneq \label{riemann-tensor-contraction}
R_\varphi . \xi  =  R^i_{k \ell j} \, \varphi^k \, \varphi^\ell \, \xi_i \, e^j 
\monend
and $ \, (R_\varphi . \xi)_j =  R^i_{k \ell j} \, \varphi^k \, \varphi^\ell \, \xi_i $. 

\smallskip \noindent
The time derivative  
of a state $ \, q(t) \, $ on the manifold defines a contravariant vector field $ \, \zeta \, $ according to
\moneq \label{zeta-dqsurdt} 
\zeta = {{\dd q}\over{\dd t}} =  \Big( {{\dd}\over{\dd t}} q^j  \Big) \, e_j  \equiv  \dot {q}^j \, e_j  
\monend 

\smallskip \noindent 
and $ \, \zeta^j =  \dot {q}^j $. 
In a similar way, the first order time derivative of a covector $\, \xi = \xi_j \, e^j  \, $
along a trajectory $ \, q(t) \, $  satisfies the conditions 
$ \, {{\dd \xi}\over{\dd t}} = \big(  \dot {\xi_j} -  \Gamma^k_{j \ell} \, \xi_k \, \zeta^\ell   \big) \, e^j $.

\bigskip  \noindent

{\bf Proposition 1.} 
  {\bf Variation of the first and second order time derivatives of a state on a Riemannian manifold}

\smallskip \noindent
We consider a given trajectory position $ \, q(t) \, $ on a Riemannian manifold $ \, Q $.
We denote  the velocity tangent vector by $ \, \zeta = {{\dd q}\over{\dd t}} $.
This trajectory position is supposed
to vary in an infinitesimal way with the variation $ \, \delta q = \delta q^j \, e_j \, $ 
of the state. We have the relations 
\moneq \label{delta-dt-q} 
\delta \Big( {{\dd q}\over{\dd t}} \Big) =  \delta \zeta =
\Big[ \delta(\zeta^j) + \Gamma^j_{k \ell} \, \zeta^\ell \, \delta q^k  \Big] \, e_j 
\monend 
\moneq \label{delta-d2t-q} 
\delta \Big( {{\dd^2 q}\over{\dd t^2}} \Big) =  \delta \Big( {{\dd \zeta}\over{\dd t}} \Big) =
\Big[ \delta ( {\dot \zeta}^j )
+ 2 \,    \Gamma^j_{k \ell}  \, \zeta^k \, \delta(\zeta^\ell)
+ \Big(\partial_k \Gamma^j_{\ell m}  \, \zeta^\ell \, \zeta^m
+ \Gamma^j_{k \ell}  \, \Big( {{\dd \zeta}\over{\dd t}} \Big)^\ell \Big) \, \delta q^k 
\Big] \, e_j
\monend 
%

\smallskip \monitem Proof of Proposition 1 

\noindent
The relation (\ref{delta-dt-q}) is an easy consequence of the variation
$ \, \delta e_\ell = \Gamma^j_{k \ell} \, \delta q^k \, e_j \, $
of a tangent vector in some infinitesimal variation.
We have also $ \, {{\dd \zeta}\over{\dd t}} = 
[ \dot{\zeta^j}  + \Gamma^j_{k \ell} \, \zeta^\ell \, \zeta^k  ] \, e_j $. Then we have 

\smallskip \noindent
$ \delta \big(  {{\dd \zeta}\over{\dd t}} \big) =
\big[ \delta  ( {\dot \zeta}^j ) +  (\partial_k \Gamma^j_{m \ell} ) \, \delta q^k \, \zeta^\ell \, \zeta^m
  +  2 \,\,  \Gamma^j_{k \ell} \, \zeta^m \, \delta (\zeta^\ell) \big] \, e_j
+ \big( {{\dd \zeta}\over{\dd t}} \big)^\ell  \, \Gamma^j_{k \ell} \, \,  \delta q^k \, e_j $

\smallskip \noindent
and the relation (\ref{delta-d2t-q})  is established. \hfill $\square$ 

\bigskip \noindent 
{\bf Proposition 2.}
{\bf Second time derivative of a covariant vector}

\smallskip \noindent
If  $\, \xi = \xi_j \, e^j \, $ is a covector field on a manifold $ \, Q $, 
we can explicit the components  of the  second time derivative
$ \, {{\dd^2 \xi}\over{\dd t^2}} = \big(  {{\dd^2 \xi}\over{\dd t^2}} \big)_j \, e^j \, $ 
of this co-vector along a trajectory position $ \, q(t) $: 
\moneq \label{d2xi-sur-dt2} 
 {{\dd^2 \xi}\over{\dd t^2}} = \Big[  \ddot {\xi_j}
- 2 \,  \Gamma^k_{j \ell} \, \Big({{\dd \xi}\over{\dd t}}\Big) _k   \, \zeta^l
- \Gamma^k_{j \ell} \, \xi_k \, \Big({{\dd \zeta}\over{\dd t}}\Big)^\ell
+ \Big( R^k_{\ell m j} - \partial_j  \Gamma^k_{\ell m} \Big) \, \xi_k \, \zeta^\ell \, \zeta^m \Big] \, e^j  \, .
\monend 
%

\smallskip \monitem Proof of Proposition 2. 

\noindent   
We differentiate relatively to time the first order derivative
$  {{\dd \xi}\over{\dd t}} \!=\! \big(  \dot {\xi_j} \!-\!  \Gamma^k_{j \ell} \, \xi_k \, \zeta^\ell   \big) \, e^j $. 
Then

\noindent
$   {{\dd^2 \xi}\over{\dd t^2}} =  {{\dd}\over{\dd t}} \big(  \dot {\xi_j} -  \Gamma^k_{j \ell} \, \xi_k \, \zeta^\ell   \big) \, e^j
+ \big(   {{\dd \xi}\over{\dd t}} \big)_k  \,  {{\dd e^k}\over{\dd t}} $

\smallskip  \quad $ \,\,\, =
\Big[  \ddot {\xi_j} - \big( \partial_m   \Gamma^k_{j \ell} \big) \,  \xi_k \, \zeta^\ell \, \zeta^m
-  \Gamma^k_{j \ell} \,  \dot {\xi_k} \, \zeta^\ell
-  \Gamma^k_{j \ell} \,    \xi_k \,  \dot { \zeta^\ell}   - \big(   {{\dd \xi}\over{\dd t}} \big)_k \,  \Gamma^k_{j \ell}  \, \zeta^\ell \Big] \, e^j  $ 

\smallskip \quad $ \,\,\, =
\Big[ \ddot {\xi_j} - \big( \partial_m   \Gamma^k_{j \ell} \big) \,  \xi_k \, \zeta^\ell \, \zeta^m
   -  \Gamma^k_{j \ell} \, \Big( \big( {{\dd \xi}\over{\dd t}} \big)_k +   \Gamma^m_{k p}  \, \xi_m \,\, \zeta^p \Big) \, \zeta^\ell 
 -  \Gamma^k_{j \ell} \,    \xi_k \, \Big( \big( {{\dd \zeta}\over{\dd t}} \big)^\ell -  \Gamma^\ell_{p q}  \, \zeta^p \, \zeta^q \Big) $

\qquad \qquad    $   -  \Gamma^k_{j \ell}  \,  \big( {{\dd \xi}\over{\dd t}} \big)_k \,  \zeta^\ell \Big] \, e^j $

\smallskip \quad $ \,\,\, =
\Big[ \ddot {\xi_j} - 2 \, \Gamma^k_{j \ell}  \,  \big( {{\dd \xi}\over{\dd t}} \big)_k \,  \zeta^\ell
  -  \Gamma^k_{j \ell} \,    \xi_k \,  \big( {{\dd \zeta}\over{\dd t}} \big)^\ell $
     $ + \big( -   \partial_m   \Gamma^k_{j \ell} + \Gamma^s_{m \ell}  \, \Gamma^k_{s j}
-    \Gamma^s_{j \ell} \, \Gamma^k_{s m} \big) \, \xi_k \, \zeta^\ell \, \zeta^m \Big] \, e^j $. 

\smallskip \noindent 
But thanks to (\ref{riemann-tensor}), we have 
$ \,\, R^k_{\ell m j} =   \partial_j  \Gamma^k_{\ell m}  -   \partial_m   \Gamma^k_{j \ell}
+ \Gamma^s_{m \ell}  \, \Gamma^k_{s j} -    \Gamma^s_{j \ell} \, \Gamma^k_{s m}   \,\, $

\smallskip \noindent
and we deduce that $ \,\,  -   \partial_m   \Gamma^k_{j \ell} + \Gamma^s_{m \ell}  \, \Gamma^k_{s j} -    \Gamma^s_{j \ell} \, \Gamma^k_{s m}  = 
R^k_{\ell m j} - \partial_j  \Gamma^k_{\ell m} $.

\smallskip  \noindent 
Then $  \,\, {{\dd^2 \xi}\over{\dd t^2}} = \big[  \ddot {\xi_j}
- 2 \,  \Gamma^k_{j \ell} \, \big({{\dd \xi}\over{\dd t}}\big) _k   \, \zeta^l
- \Gamma^k_{j \ell} \, \xi_k \, \big({{\dd \zeta}\over{\dd t}}\big)^\ell
+ \big( R^k_{\ell m j} - \partial_j  \Gamma^k_{\ell m} \big) \, \xi_k \, \zeta^\ell \, \zeta^m \big] \, e^j \,\, $ 

\smallskip  \noindent 
and the property is established. \hfill $\square$

\bigskip

\bigskip    \bigskip  \noindent {\bf \large    3) \quad  Pontryagin  framework for differential equations}    

\smallskip \noindent
We consider a dynamical system in a finite dimensional euclidian space.
A state vector~$ \, y(t ,\, \lambda) \in  \R^d  \,$
is submitted to a system of  first order differential equations
\moneq \label{edo} 
{{{\rm d}y}\over{{\rm d}t}} =  f(y(t) ,\, \lambda (t) ,\, t) .
\monend      
This system is controlled by a set of dynamical parameters   $\, \lambda (t) $. 
The  initial condition   takes the form 
$ \, y(0 ,\, \lambda)  =  x $. 
We search an optimal solution that minimizes  the  cost function 
\moneq \label{cout-pontryaguine} 
J( \lambda) \, \equiv \, 
\int_{0}^{T} \, g \big( y(t),\, \lambda(t),\,t \big) \, {\rm d}t
\monend
Pontryagin's main idea (see {\it e.g.} \cite{PBGM62})
can be formulated as follows. 
Consider the differential equation
$ \, {{{\rm d}y}\over{{\rm d}t}} =  f(y(t) ,\, \lambda (t) ,\, t) \, $ 
as a  constraint satisfied by the variable  $\, y \, $ and
introduce a Lagrange multiplier  $ \,   p = p(t) \, $ 
associated with this constraint.
Then a Lagrangian functional
\moneqstar
{\cal L}(y ,\,\lambda ,\,p) \equiv 
\int_{0}^{T} \, g(y,\,\lambda,\,t) \, {\rm d}t \,+\, \int_{0}^{T}  p(t) \, \Bigl(
{{{\rm d}y}\over{{\rm d}t}} - f(y ,\, \lambda ,\, t) \Bigr) \, {\rm d}t 
\monendstar
is naturally associated with the cost function and the differential equation
viewed as a constraint.
After a classical integration by parts of the variation $ \, \delta  {\cal L} \, $ of the Lagrangian
(see {\it e.g.} \cite {PBGM62}), it is well known that if the adjoint state   $\, p(t) \,$ satisfies
the  following adjoint equation
\moneqstar
{{{\rm d}p}\over{{\rm d}t}} +  p \, {{\partial f}\over{\partial y}}   
- {{\partial g}\over{\partial y}}  =  0 
\monendstar
and the    final condition: $ \, p \,(T) =  0 $, 
then the variation $\,  \delta J  \, $ of the cost function is given by the relation 
\moneqstar
\delta J =  \int_{0}^{T} \, \Bigl[ \,  {{\partial g}\over{\partial  \lambda}}
- p \, {{\partial f}\over{\partial \lambda}}  \, \Bigl] \,  \delta \lambda (t)   \,  {\rm d}t   
\monendstar
for a given variation $\, \delta \lambda  \, $
of the parameter. 
At the optimum this variation is identically null and this is expressed with the 
Pontryagin optimality condition
$ \, {{\partial g}\over{\partial  \lambda}} - p \, {{\partial f}\over{\partial \lambda}}  = 0 $.

\bigskip   \bigskip   \noindent {\bf \large    4) \quad  Optimal dynamics for a quadratic cost functional}    

\smallskip \noindent
We consider now a mechanical system described by a state $ \, q(t) \, $ on a  manifold $ \, Q \, $
of finite dimension. We suppose given a mechanical Lagrangian 
\moneqstar
L(q, \, \zeta ) \equiv   K(q, \zeta) - V(q)
\monendstar
with $ \, K(q,\, \zeta) \equiv M_{k \ell}(q) \, \zeta^k \, \zeta^\ell \, $ the kinetic energy
of the system.
It defines a metric through the mass
matrix as observed previously in (\ref{energie-cinetique}).
The Euler-Lagrange equations of a free evolution take the form 
\moneqstar
{{\rm d}\over{{\rm d}t}} \Big( {{\partial L}\over{\partial \zeta^i}} \Big) 
=  {{\partial L}\over{\partial q^i}} 
\monendstar
for all degrees of freedom.
These equations take a Riemannian form:
\moneq \label{euler-lagrange} 
M_{k \ell} \, \big( {\dot \zeta}^\ell + \Gamma^\ell_{i j} \, \zeta^i \, \zeta^j \big) + \partial_k V = 0  
\monend
and the proof of this relation can be found in \cite{DFRV15,Juan-these-2013}. 
After some index juggling,   
the relation (\ref{euler-lagrange}) can be written 
$ \,\,  {\dot \zeta}^j + \Gamma^j_{k \ell} \,  \zeta^k \, \zeta^\ell + M^{j \ell} \,   \partial_\ell V = 0 $. 

\smallskip \noindent 
The objective of an engineering process is the control of the state $ \, q(t) \, $ along the time,
adding forces and torques $ \,\, u = u_k \, e^k \,\, $ to the natural evolution.
Observe that the control source $ \, u \, $ is a covariant vector field on the manifold. 
We obtain with this process (see {\it e.g.} \cite{Juan-these-2013})
the evolution equations
\moneqstar
M_{k \ell} \, \big( {\dot \zeta}^\ell + \Gamma^\ell_{i j} \, \zeta^i \, \zeta^j \big) + \partial_k V = u_k . 
\monendstar 
We can introduce the contravariant components $ \, u^j = M^{j k} \, u_k \, $ for the covector.
Then the dynamical  evolution equations can be written as 
\moneq \label{euler-lagrange-avec-forces}
{\dot \zeta}^j + \Gamma^j_{k \ell} \,  \zeta^k \, \zeta^\ell + M^{j \ell} \,   \partial_\ell V = u^j \, . 
\monend
A fundamental idea of our approach \cite{LV95} is to enforce the coherence 
of the controlled mechanical system with
a cost function $ \, J(u) \, $ that respects the Riemannian structure of the free evolution. 
The choice of a quadratic functional is proposed in \cite{Juan-these-2013}: 
\moneq \label{cout-quadratique}
J(u) =  {1\over2} \, \int_0^T M_{k \ell} (q) \, u^k \, u^\ell \, {\rm d}t
\monend 
It is possible to make a link with the Pontryagin approach (\ref{edo})(\ref{cout-pontryaguine})
with the choice proposed in \cite{DFRV15}:
\moneqstar
y =  \{ q^j \,, \, \zeta^j  \} \,,\,\,  
f =   \{ \zeta^j \,, \, 
- \Gamma^j_{k \ell} \,  \zeta^k \, \zeta^\ell - M^{j \ell} \, \partial_\ell V + u^j \} \,,\,\,     
\lambda  =  \{ u^k \} \,,\,\,    
 g =  {1\over2} \, M_{k \ell} (q) \, u^k \, u^\ell .
\monendstar
Observe that the quadratic functional (\ref{cout-quadratique}) has an intrinsic 
structure that respects the fundamental mechanical constraints. 
The Lagrange multiplers or adjoint states take the form
$\,\, p =  \{ \rho_j \,, \, \xi_j \} \, $
with $ \, \rho = \rho_j \, e^j \, $ associated with the first equation
$ \, {{\dd q}\over{\dd t}} = \zeta \, $ and
$ \, \xi = \xi_j \, e^j \, $ multiplying  the dynamics
$ \,\,  {\dot \zeta}^j + \Gamma^j_{k \ell} \,  \zeta^k \, \zeta^\ell + M^{j \ell} \,   \partial_\ell V - u^j = 0 $. 
A very beautiful result established in \cite{Juan-these-2013} is the interpretation of the adjoint state 
$ \, \xi \, $ as exactly equal to the forces and torques.
We have  
\moneqstar
\xi = u ,
\monendstar
{\it id est}
$ \, \xi_k = u_k \, $ for all the covariant components. 
Moreover, a precise evolution equation for the dual variable has been established.

\bigskip  \noindent

{\bf Theorem 1.}
{\bf Covariant evolution equation  of the optimal force}

\smallskip \noindent
With the above notations and hypotheses, the forces and torques  $ \, u \,$ 
satisfy the following time evolution:
\moneq \label{evolution-xi-old} 
\Big( {{{\rm d}^2 u}\over{{\rm d}t^2}} \Big)_j +
R^i_{k \ell j} \,   {\dot q}^k \,{\dot q}^\ell \, u_i  + 
\big( \nabla^2_{j k} V \big) \, u^k = 0 \, . 
\monend 

\smallskip \noindent
This relation has been derived in Rojas-Quintero's  thesis \cite{Juan-these-2013},
and is  presented in 
\cite{DFRV15}.  

\smallskip \smallskip \noindent
One fundamental case is the double pendulum and it has been considered 
for an experimental confrontation.
In this  case, the manifold $ \, Q \, $ is of dimension 2. 
The efficiency of the choice of a covariant quadradic functional is not
{\it a priori} obvious. It is studied for the double pendulum
and compared with experiments and simulations in the references 
\cite{RQ-VC-S-21} and \cite{RDR22}. 

\bigskip   \bigskip   \noindent {\bf \large    5) \quad  General second order covariant adjoint equation}    

\smallskip \noindent
We consider in this contribution a general cost function
\moneq \label{cost-function} 
J(u) = \int_0^T \gamma( q,\, \zeta ,\, u) \, \dd t
\monend
instead of the quadratic functional (\ref{cout-quadratique}).
The Lagrangian of the problem introduces the adjoint states $ \, \rho \, $ and $ \, \xi \, $
relative to each equation of the dynamical system
\moneq \label{dynamical-system} 
  {{\dd q}\over{\dd t}} = \zeta \,,\,\,   {{\dd \zeta}\over{\dd t}} - \psi(q) =  u 
\monend
and we have
\moneq \label{lagrangien} 
{\cal L} = J(u) + \int_0^T \rho \, \Big( {{\dd q}\over{\dd t}} - \zeta \Big)\, \dd t
+ \int_0^T \xi \, \Big(   {{\dd \zeta}\over{\dd t}} - \psi(q) - u \Big)\, \dd t  .
\monend  

\bigskip \noindent 
\newpage
\noindent 
{\bf Proposition  3.}
{\bf Variation of the Lagrangian}

\noindent
For arbitrary variations  $ \, (\delta q,\,  \delta \zeta )  \, $  of the state
  $ \, (q ,\, \zeta) $,  
$ \, (\delta p ,\, \delta \xi )  $ of the Lagrange multipliers~$ \, p \, $ and $ \, \xi $,
and $ \, \delta u \, $ of the control  variable $ \, u $, 
we have the following variation $ \, \delta  {\cal L} \, $ of the lagrangian defined in~(\ref{lagrangien}): 
\moneq \label{delta-Lagrangien} 
\delta  {\cal L}  = \left\{\begin{array} {l} \displaystyle
\int_0^T \delta \rho \, \Big( {{\dd q}\over{\dd t}} - \zeta \Big)  \dd t
+ \int_0^T \delta \xi \, \Big(   {{\dd \zeta}\over{\dd t}} - \psi(q) - u \Big)  \dd t \\  \vspace{-3 mm} \\ \displaystyle
+ \Big[ \Big(  {{\partial \gamma}\over{\partial \zeta}} - {{\dd \xi}\over{\dd t}} \Big) \, \delta q 
  + \xi \,  \delta \zeta \Big]_0^T
+ \int_0^T \Big( {{\partial \gamma}\over{\partial u}} - \xi \Big) \, ( \delta u )^j \,  \dd t \\ \vspace{-3 mm} \\ \displaystyle 
+\int_0^T \Big[ \Big( {{\partial \gamma}\over{\partial q}} \Big)_j
- \Big(  {{\dd}\over{\dd t}} \Big(  {{\partial \gamma}\over{\partial \zeta}} \Big) \Big)_j 
+ \Big( {{\dd^2 \xi}\over{\dd t^2}} \Big)_j 
\\ \vspace{-3 mm} \\ \displaystyle \qquad  \qquad 
+ R^k_{\ell j m} \,  \xi_k \, \zeta^\ell \, \zeta^m  
- \xi_k \, (\partial_j \psi )^k \Big] \, \delta q^j \, \dd t  \,,
\end{array} \right. \monend
where $ \,  R^k_{\ell j m} \, $ is the Riemann curvature tensor defined in (\ref{riemann-tensor}).

\smallskip \monitem Proof of  Proposition 3. 

\noindent
The Lagrangien of this problem can be written
$ \,   {\cal L} = J(u) +  {\cal L}_1 +  {\cal L}_2 \, $ with
\moneq \label{L1-et-L2} 
{\cal L}_1 = \int_0^T \rho \, \Big( {{\dd q}\over{\dd t}} - \zeta \Big)\, \dd t \,,\,\,
{\cal L}_2 = \int_0^T \xi \, \Big(   {{\dd \zeta}\over{\dd t}} - \psi(q) - u \Big)\, \dd t .
\monend
Recall that we have $ \, \rho = \rho_j \, e^j $, $ \, \xi = \xi_j \, e^j $,
$\, (\delta q)^j = \delta(q^j) $, 
$ \, (\delta \zeta)^j = \delta(\zeta^j) +  \Gamma^j_{k \ell} \,\, \zeta^k \, \delta q^\ell \, $  
and $ \, (\delta u)^j = \delta(u^j) +  \Gamma^j_{k \ell} \,\, u^k \, \delta q^\ell $. 
We take the variation of the three terms of the Lagrangian function. For the cost function
defined in (\ref{cost-function}), we have

\smallskip \noindent 
$ \delta J =  \int_0^T \big[ {{\partial \gamma}\over{\partial q}} \, \delta q
+ {{\partial \gamma}\over{\partial \zeta}} \, \delta \zeta
+ {{\partial \gamma}\over{\partial u}} \, \delta u \big] \, \dd t $

 \smallskip \quad 
$ \,\, =  \int_0^T \big[ \big({{\partial \gamma}\over{\partial q}} \big)_{\!j} \,  \delta q^j
+ \big({{\partial \gamma}\over{\partial \zeta}} \big)_{\!j}  \,  (\delta \zeta)^j   
+ \big({{\partial \gamma}\over{\partial u}} \big)_{\!j}  \,  (\delta u)^j   \big] \, \dd t $

 \smallskip \quad 
$ \,\, =  \int_0^T \big[ \big({{\partial \gamma}\over{\partial q}} \big)_{\!j} \,  \delta q^j
+ \big({{\partial \gamma}\over{\partial \zeta}} \big)_{\!j}  \,
\big(  \delta(\zeta^j) +  \Gamma^j_{k \ell} \,\, \zeta^k \, \delta q^\ell \big) 
+ \big({{\partial \gamma}\over{\partial u}} \big)_{\!j}  \,  (\delta u)^j   \big] \, \dd t $

\smallskip \noindent and 
\moneq \label{delta-J} 
\delta J =   \int_0^T \Big[ \Big({{\partial \gamma}\over{\partial q}} \Big)_{\!j}
  + \Gamma^\ell_{k j} \,\,  \Big({{\partial \gamma}\over{\partial \zeta}} \Big)_{\!\ell} \, \zeta^k  \Big] \,  \delta q^j \, \dd t
+  \int_0^T \Big[  \Big({{\partial \gamma}\over{\partial \zeta}} \Big)_{\!j} \,  \delta(\zeta^j)
+ \Big({{\partial \gamma}\over{\partial u}} \Big)_{\!j}  \,  (\delta u)^j   \Big] \, \dd t  .
\monend
From $\, {{\dd q}\over{\dd t}} = \dot{q}^j \, e_j = \zeta^j \, e_j $, we have

\smallskip \noindent 
$   \delta \big({{\dd q}\over{\dd t}}\big) = \delta \dot{q}^j \, e_j + \dot{q}^j \, \delta e_j $
$ = \big(  \delta \dot{q}^j + \Gamma^j_{k \ell} \,\, \zeta^k \, \delta q^\ell  \big) \, $
$ = \big(  \delta (\zeta^j)  + \Gamma^j_{k \ell} \,\, \zeta^k \, \delta q^\ell  \big) \, $
$ =  ( \delta \zeta )^j  \, $

\smallskip \noindent 
and by recalling (\ref{delta-dt-q}) of  Proposition~1, 

\smallskip \noindent
$ \,  \delta \big( \rho \, ( {{\dd q}\over{\dd t}} - \zeta ) \big) =  \delta \rho \, \big( {{\dd q}\over{\dd t}} - \zeta \big)
+ \rho_j \, \big(  \delta \dot{q}^j + \Gamma^j_{k \ell} \,\, \zeta^k \, \delta q^\ell  -   ( \delta \zeta )^j   \big) $

\smallskip \noindent \qquad \qquad \qquad 
$ \, =  \delta \rho \, \big( {{\dd q}\over{\dd t}} - \zeta \big) + {{\dd}\over{\dd t}} \big(  \rho_j \,   \delta q^j \big)
- \dot{\rho}_j  \, \delta q^j - \rho_j \,  \delta (\zeta^j) $. 

\smallskip \noindent
Then integrating by parts 

\smallskip \noindent 
$ \delta {\cal L}_1 =   \int_0^T \delta \rho \, \big( {{\dd q}\over{\dd t}} - \zeta \big)\, \dd t
+ [ \rho_j \, \delta q^j ]_0^T    - \int_0^T \dot{\rho}_j  \, \delta q^j \, \dd t 
- \int_0^T \rho^j  \, \delta (\zeta^j) \, \dd t \,\, $

\smallskip \noindent 
and 
\moneq \label{delta-L1} 
\delta  {\cal L}_1 = \big[ \rho_j \, \delta q^j \big]_0^T +
\int_0^T \delta \rho \, \Big( {{\dd q}\over{\dd t}} - \zeta \Big)\, \dd t
-  \int_0^T  \Big( \dot{\rho}_j  \, \delta q^j  + \rho_j  \, \delta (\zeta^j) \Big) \, \dd t \, . 
\monend

\smallskip \noindent 
We observe now that we have for the contravariant vector field
$ \, \delta \psi = \big( \delta \psi^j +  \Gamma^j_{k \ell} \,\, \psi^k \, \delta q^\ell \big) \, e_j $. 
We keep the compact expression $ \, \delta u =  \big( \delta u \big)^j \, e_j $. We can develop
the third term: 

\smallskip \noindent 
$ \delta {\cal L}_2 = \int_0^T \delta \xi \, \big(   {{\dd \zeta}\over{\dd t}} - \psi(q) - u \big)\, \dd t
+ \int_0^T \xi \,  \big(  \delta  {{\dd \zeta}\over{\dd t}} -    \delta \psi(q) -   \delta u \big)\, \dd t $ 

\smallskip \noindent
and from (\ref{delta-d2t-q}) and Lemma 2, we have 

\smallskip \noindent 
$ \delta {\cal L}_2 = \int_0^T \delta \xi \, \big(   {{\dd \zeta}\over{\dd t}} - \psi(q) - u \big)\, \dd t 
+  \int_0^T \xi_j \, \big[  \delta  {\dot \zeta}^j 
+ \big( \partial_k \Gamma^j_{\ell m} \, \zeta^\ell \, \zeta^m + \Gamma^j_{k \ell} \, \big( {{\dd \zeta}\over{\dd t}} \big)^\ell \big) \, \delta q^k
+ 2 \, \Gamma^j_{k \ell}  \, \zeta^k \, \delta (\zeta^\ell) \big] \,  \dd t $

\smallskip \noindent \qquad \quad 
$ -  \int_0^T \xi_j \,  \big( \partial_\ell \psi^j + \Gamma^j_{k \ell} \, \psi^k \big) \, \delta q^\ell \, \dd t 
 -  \int_0^T \xi_j \, \big( \delta u \big)^j \,  \dd t \,\,\, $ and 
\moneq \label{delta-L2} 
\!\!\!\!\! \delta  {\cal L}_2 = \left\{\begin{array} {l} 
\int_0^T \delta \xi \, \big(   {{\dd \zeta}\over{\dd t}} - \psi(q) - u \big)\, \dd t
+ \big[ \xi_j \, \delta (\zeta^j) \big]_0^T 
+ \int_0^T \big( -  {\dot \xi}_j + 2 \, \Gamma^k_{j \ell}  \, \xi_k \, \zeta^\ell \big) \, \delta (\zeta^\ell) \,  \dd t \\
\!\! +  \int_0^T \big[ \xi_k \, \big( \partial_j \Gamma^k_{\ell m} \, \zeta^\ell \, \zeta^m
  +  \Gamma^k_{j \ell} \, \big( {{\dd \zeta}\over{\dd t}} \big)^\ell - \partial_j \psi^k -  \Gamma^k_{j \ell} \,\psi^\ell 
  \big)  \big] \, \delta q^j \,  \dd t -  \int_0^T \xi_j \, \big( \delta u \big)^j   \dd t .
\end{array} \right. 
\monend
We can now add the three contributions detailed in the relations (\ref{delta-J}), (\ref{delta-L1}) and (\ref{delta-L2}):

\smallskip \noindent 
$ \delta {\cal L} =  [ \rho_j \, \delta q^j + \xi_j \, \delta (\zeta^j) ]_0^T + \int_0^T \delta \rho \, \big( {{\dd q}\over{\dd t}} - \zeta \big) \, \dd t
+ \int_0^T \delta \xi \, \big(   {{\dd \zeta}\over{\dd t}} - \psi(q) - u \big) \, \dd t $ 

\smallskip  \quad $\,\,\,$
$ + \int_0^T \big[ \big( {{\partial \gamma}\over{\partial q}} \big)_j
+  \Gamma^\ell_{k j} \, \big( {{\partial \gamma}\over{\partial \zeta}} \big)_\ell \,\, \zeta^k - {\dot \rho}_j
+ (\partial_j \Gamma^k_{\ell m}) \, \xi_k \, \zeta^\ell \, \zeta^m 
+ \Gamma^k_{j \ell} \, \xi_k \, \big(   {{\dd \zeta}\over{\dd t}} \big)^\ell - \xi_k \, \big(  \partial_j \psi^k
+  \Gamma^k_{j \ell} \,\psi^\ell \big) \big] \, \delta q^j  \, \dd t $ 

\smallskip  \quad $\,\,\,$
$ +  \int_0^T \big[ \big( {{\partial \gamma}\over{\partial u}} \big)_j - \xi_j \big] \, \big( \delta u \big)^j  \dd t
+  \int_0^T \big[ \big( {{\partial \gamma}\over{\partial \zeta}} \big)_j - \rho_j -  {\dot \xi}_j + 2 \,  \Gamma^k_{j \ell} \,\xi_k \, \zeta^\ell \big]
\, \delta ( \zeta^j ) \, \dd t $. 

\smallskip \noindent
Because $ \, \delta ( \zeta^j ) = \delta {\dot q}^j = {{\dd}\over{\dd t}} (\delta q^j) $, we can integrate by parts the last term and we obtain 

\smallskip \noindent 
$ \delta {\cal L} =  [ \rho_j \, \delta q^j + \xi_j \, \delta (\zeta^j) ]_0^T + \int_0^T \delta \rho \, \big( {{\dd q}\over{\dd t}} - \zeta \big) \, \dd t
+ \int_0^T \delta \xi \, \big(   {{\dd \zeta}\over{\dd t}} - \psi(q) - u \big) \, \dd t $ 

\smallskip  \quad $\,\,\,$
$ + \int_0^T \big[ \big( {{\partial \gamma}\over{\partial q}} \big)_j
  +  \Gamma^\ell_{k j} \, \big( {{\partial \gamma}\over{\partial \zeta}} \big)_\ell \,\, \zeta^k
  - {\dot \rho}_j + (\partial_j \Gamma^k_{\ell m}) \, \xi_k \, \zeta^\ell \, \zeta^m 
+ \Gamma^k_{j \ell} \, \xi_k \, \big(   {{\dd \zeta}\over{\dd t}} \big)^\ell - \xi_k \, \big(  \partial_j \psi^k
+  \Gamma^k_{j \ell} \,\psi^\ell \big) \big] \, \delta q^j  \, \dd t $ 

\smallskip  \quad $\,\,\,$
$ + \int_0^T \big[  \big( {{\partial \gamma}\over{\partial u}} \big)_j - \xi_j \big] \, (\delta u)^j \, \dd t $ 
$ + \big[ \big( \big( {{\partial \gamma}\over{\partial \zeta}} \big)_j - \rho_j - {\dot \xi}_j +  2 \,  \Gamma^k_{j \ell} \,\xi_k \, \zeta^\ell \big) 
  \, \delta q^j \big]_0^T $ 

\smallskip  \quad $\,\,\,$
$- \int_0^T  {{\dd}\over{\dd t}}  \big[ \big( {{\partial \gamma}\over{\partial \zeta}} \big)_j - \rho_j - {\dot \xi}_j
+  2 \,  \Gamma^k_{j \ell} \,\xi_k \, \zeta^\ell \big]  \, \delta q^j  \, \dd t  $

\smallskip  \quad 
$ =  [ \big( \big( {{\partial \gamma}\over{\partial \zeta}} \big)_j - {\dot \xi}_j +  2 \,  \Gamma^k_{j \ell} \,\xi_k \, \zeta^\ell \big) \big) \, \delta q^j 
  + \xi_j \, \delta (\zeta^j) ]_0^T + \int_0^T \delta \rho \, \big( {{\dd q}\over{\dd t}} - \zeta \big) \, \dd t
+ \int_0^T \delta \xi \, \big(   {{\dd \zeta}\over{\dd t}} - \psi(q) - u \big) \, \dd t $

\smallskip  \quad $\,\,\,$
$ +  \int_0^T \big[  \big( {{\partial \gamma}\over{\partial u}} \big)_j - \xi_j \big] \, (\delta u)^j \, \dd t
+  \int_0^T \big[ \big( {{\partial \gamma}\over{\partial q}} \big)_j
  +  \Gamma^\ell_{k j} \, \big( {{\partial \gamma}\over{\partial \zeta}} \big)_\ell \, \zeta^k 
  -  {{\dd}\over{\dd t}}  \big( {{\partial \gamma}\over{\partial \zeta}} \big)_j + \ddot{\xi}_j
  +  (\partial_j \Gamma^k_{\ell m}) \, \xi_k \, \zeta^\ell \, \zeta^m $

\smallskip  \quad $\,\,\,$
$   +  \, \Gamma^k_{j \ell} \, \xi_k \, \big( {{\dd \zeta}\over{\dd t}}\big)^\ell
- 2 \, {{\dd}\over{\dd t}} \big( \Gamma^k_{j \ell} \, \xi_k \, \zeta^\ell \big) 
- \xi_k \, \big(  \partial_j \psi^k +  \Gamma^k_{j \ell} \,\psi^\ell \big) \big] \, \delta q^j  \, \dd t $. 

\smallskip \smallskip \noindent
The boundary term can be simplified:

\smallskip \noindent 
$ \,  \big[ \big( \big( {{\partial \gamma}\over{\partial \zeta}} \big)_j - \dot{\xi}_j  
  +  2 \,  \Gamma^k_{j \ell} \,\xi_k \, \zeta^\ell \big) \big) \, \delta q^j   + \xi_j \, \delta (\zeta^j) \big]_0^T =
    \big[ \big( \big( {{\partial \gamma}\over{\partial \zeta}} \big)_j - \big( {{\dd \xi}\over{\dd t}}\big)_j \big) \, \delta q^j
      + \xi_j \, \big(  \delta (\zeta^j) +  \Gamma^j_{k \ell} \, \zeta^k \,  \delta q^\ell \big)  \big]_0^T  $

\smallskip \qquad \qquad 
$ =    \big[ \big( \big( {{\partial \gamma}\over{\partial \zeta}} \big)_j - \big( {{\dd \xi}\over{\dd t}} \big)_j \big) \, \delta q^j
+  \xi_j \,  (  \delta \zeta ) ^j \big]_0^T
=  \big[ \big({{\partial \gamma}\over{\partial \zeta}} - {{\dd \xi}\over{\dd t}} \big) \, \delta q 
  + \xi \,  \delta \zeta \big]_0^T $

\smallskip \noindent
and natural  boundary conditions 
are put in evidence.

\smallskip \noindent 
We focus now our attention on the term containing the variation
$ \, \delta q \, $ in factor. We have

\smallskip  \noindent 
$ \int_0^T \big[ \big( {{\partial \gamma}\over{\partial q}} \big)_j
  +  \Gamma^\ell_{k j} \, \big( {{\partial \gamma}\over{\partial \zeta}} \big)_\ell \, \zeta^k
- {{\dd}\over{\dd t}} \big(  {{\partial \gamma}\over{\partial \zeta}} \big)_{\!j} 
  + \ddot{\xi}_j
  +  (\partial_j \Gamma^k_{\ell m}) \, \xi_k \, \zeta^\ell \, \zeta^m  
  +  \, \Gamma^k_{j \ell} \, \xi_k \, \big( {{\dd \zeta}\over{\dd t}}\big)^\ell $

\smallskip  \quad $\,\,\,$
$ - 2 \, {{\dd}\over{\dd t}} \big( \Gamma^k_{j \ell} \, \xi_k \, \zeta^\ell \big) 
- \, \xi_k \, \big(  \partial_j \psi^k +  \Gamma^k_{j \ell} \,\psi^\ell \big) \big] \, \delta q^j  \, \dd t  $

\smallskip  \quad 
$ =  \int_0^T \big[  \big( {{\partial \gamma}\over{\partial q}} \big)_j
- \big(  {{\dd}\over{\dd t}} \big(  {{\partial \gamma}\over{\partial \zeta}} \big) \big)_j 
+ {\ddot \xi}_j
+    (\partial_j \Gamma^k_{\ell m}) \, \xi_k \, \zeta^\ell \, \zeta^m  + \Gamma^k_{j \ell} \, \xi_k \, \big( {{\dd \zeta}\over{\dd t}}\big)^\ell 
- 2 \, (\partial_m \Gamma^k_{j \ell} ) \, \xi_k \, \zeta^\ell \, \zeta^m $ 

\smallskip  \quad $\,\,\,$
$ - 2 \, \Gamma^k_{j \ell} \, \big(  \big( {{\dd \xi}\over{\dd t}} \big)_k +  \Gamma^s_{k p} \, \xi_s \, \zeta^p \big) \, \zeta^\ell
- 2 \, \Gamma^k_{j \ell} \, \xi_k \, \big(  \big( {{\dd \zeta}\over{\dd t}} \big)^\ell -  \Gamma^\ell_{s m} \, \zeta^s \, \zeta^m \big)
- \xi_k \, (\partial_j \psi )^k \big] \, \delta q^j \, \dd t  $

\smallskip  \quad 
$ =  \int_0^T \big[  \big( {{\partial \gamma}\over{\partial q}} \big)_j
- \big(  {{\dd}\over{\dd t}} \big(  {{\partial \gamma}\over{\partial \zeta}} \big) \big)_j 
+ {\ddot \xi}_j
-  \Gamma^k_{j \ell} \, \xi_k \,  \big( {{\dd \zeta}\over{\dd t}} \big)^\ell
- 2  \, \Gamma^k_{j \ell} \,  \big( {{\dd \xi}\over{\dd t}} \big)_k  \, \zeta^\ell
+ \big( \partial_j \Gamma^k_{\ell m} - 2 \, \partial_m \Gamma^k_{j \ell} \big) \,  \xi_k \, \zeta^\ell \, \zeta^m $

\smallskip  \quad $\,\,\,$
$  + \, 2 \, ( \Gamma^s_{\ell m} \, \Gamma^k_{s j} - \Gamma^s_{j \ell} \, \Gamma^k_{s m} )  \,  \xi_k \, \zeta^\ell \, \zeta^m
- \xi_k \, (\partial_j \psi )^k \big] \, \delta q^j \, \dd t  $

\smallskip  \quad
$ =  \int_0^T \big[  \big( {{\partial \gamma}\over{\partial q}} \big)_j
- \big(  {{\dd}\over{\dd t}} \big(  {{\partial \gamma}\over{\partial \zeta}} \big) \big)_j 
+ {\ddot \xi}_j
-  \Gamma^k_{j \ell} \,\, \xi_k \,  \big( {{\dd \zeta}\over{\dd t}} \big)^\ell  
- 2  \, \Gamma^k_{j \ell} \,  \big( {{\dd \xi}\over{\dd t}} \big)_k  \, \zeta^\ell
+ \big( \partial_j \Gamma^k_{\ell m} \big) \,  \xi_k \, \zeta^\ell \, \zeta^m $

\smallskip  \quad $\,\,\,$
$ + \, 2 \,  \big( R^k_{\ell m j} - \partial_j \Gamma^k_{\ell m} \big) \,  \xi_k \, \zeta^\ell \, \zeta^m  
- \xi_k \, (\partial_j \psi )^k \big] \, \delta q^j \, \dd t  $

\smallskip 
\hfill because $ \, R^k_{\ell m j} =  \partial_j \Gamma^k_{\ell m}  - \partial_m \Gamma^k_{j \ell}
+  \Gamma^s_{\ell m} \, \Gamma^k_{s j} - \Gamma^s_{j \ell} \, \Gamma^k_{s m} $

\smallskip \quad 
$ =  \int_0^T \big[ \big( {{\partial \gamma}\over{\partial q}} \big)_j
- \big(  {{\dd}\over{\dd t}} \big(  {{\partial \gamma}\over{\partial \zeta}} \big) \big)_j 
+ {\ddot \xi}_j
-  \Gamma^k_{j \ell} \,\, \xi_k \,  \big( {{\dd \zeta}\over{\dd t}} \big)^\ell  
- 2  \, \Gamma^k_{j \ell} \,  \big( {{\dd \xi}\over{\dd t}} \big)_k  \, \zeta^\ell $ 

\smallskip  \quad $\,\,\,$
$ + \, 2 \,  R^k_{\ell m j} \,  \xi_k \, \zeta^\ell \, \zeta^m -  (\partial_j \Gamma^k_{\ell m} \big) \,  \xi_k \, \zeta^\ell \, \zeta^m  
- \xi_k \, (\partial_j \psi )^k \big] \, \delta q^j \, \dd t  $

\smallskip  \quad
$ =  \int_0^T \big[ \big( {{\partial \gamma}\over{\partial q}} \big)_j 
- \big(  {{\dd}\over{\dd t}} \big(  {{\partial \gamma}\over{\partial \zeta}} \big) \big)_j 
+ \big( {{\dd^2 \xi}\over{\dd t^2}} \big)_j +  R^k_{\ell m j} \,  \xi_k \, \zeta^\ell \, \zeta^m  
- \xi_k \, (\partial_j \psi )^k \big] \, \delta q^j \, \dd t  $

\smallskip \noindent
 due to Lemma 2. We deduce a new expression for the variation of the Lagrangian:  

\smallskip \noindent 
$ \delta {\cal L} =  \big[ \big(  {{\partial \gamma}\over{\partial \zeta}} - {{\dd \xi}\over{\dd t}} \big) \, \delta q 
  + \xi \,  \delta \zeta \big]_0^T  + \int_0^T \delta \rho \, \big( {{\dd q}\over{\dd t}} - \zeta \big)  \dd t
+ \int_0^T \delta \xi \, \big(   {{\dd \zeta}\over{\dd t}} - \psi(q) - u \big)  \dd t $

\smallskip \quad $\,\,\,$  
$ + \int_0^T \big( {{\partial \gamma}\over{\partial u}}  - \xi \big) \,  \delta u  \,  \dd t
+\int_0^T \big[ \big( {{\partial \gamma}\over{\partial q}} \big)_j
- \big(  {{\dd}\over{\dd t}} \big(  {{\partial \gamma}\over{\partial \zeta}} \big) \big)_j 
+ \big( {{\dd^2 \xi}\over{\dd t^2}} \big)_j +  R^k_{\ell m j} \,  \xi_k \, \zeta^\ell \, \zeta^m  
- \xi_k \, (\partial_j \psi )^k \big] \, \delta q^j \, \dd t  $

\smallskip \noindent 
and the Proposition is established. \hfill $\square$ 

\smallskip \noindent
We observe from (\ref{delta-Lagrangien}) that the Pontryagin optimality condition
is written
\moneqstar
    {{\partial \gamma}\over{\partial u}} = \xi  .
\monendstar
The adjoint variable $ \, \xi \, $ is no more equal to the forces and torques $ \, u \, $ but the relation
between the two variables is completely explicited. 

\smallskip \noindent 
The boundary  conditions take the quite unusual form 
\moneq \label{conditions-limites} 
\Big[ \Big( {{\partial \gamma}\over{\partial \zeta}} -  {{\dd \xi}\over{\dd t}} \Big) \, \delta q
  + \xi \, \delta \zeta \Big]_0^T = 0
\monend

\smallskip \noindent 
because they can cover
several cases. To fix the ideas, when the initial conditions take the usual form
$ \, q(0) = q_0 \, $ and $ \, \zeta(0) = \zeta_0 $, with fixed given data
$\,  q_0 \, $ and $ \, \zeta_0 $, we have in consequence $ \, \delta q(0) = 0 \, $
and $ \, \delta \zeta(0) = 0 $. Then the boundary conditions (\ref{conditions-limites})
express simply a null condition at the final time: $\, \xi(T) = 0 \, $ and
$ \, \big( {{\dd \xi}\over{\dd t}} -  {{\partial \gamma}\over{\partial \zeta}} \big)(T) = 0 $. 
We can also consider for other applications that initial and final states are imposed:
$ \, q(0) = q_0 \, $ and $ \, q(T) = q_T $. In this case, 
 $ \, \delta q(0) = \delta q(T) = 0 \, $ and the expression  (\ref{conditions-limites})
express conditions for the second Lagrange multiplier at the initial and final time:
 $\, \xi(0) = \xi(T) = 0 $. Other cases can be naturally considered.

\bigskip \noindent 
{\bf Theorem 2.}
{\bf Second order adjoint evolution equation}

\noindent
When the source term derives from a potential, {\it id est} 
$ \,\,  \psi^k(q) =  -\partial^k V = \- M^{k \ell} \, \partial_\ell  V    $,
then we have no constraint for the first adjoint state $ \, \rho \, $ and 
we have a  second order dynamics for  the second Lagrange multiplier:  
\moneq \label{dynamique-second-adjoint}
{{\dd^2 \xi}\over{\dd t^2}} - R_\zeta \,.\, \xi + {{\partial \gamma}\over{\partial q}}
-  {{\dd}\over{\dd t}} \Big( {{\partial \gamma}\over{\partial \zeta}} \Big) + \nabla^2 V . \,\xi = 0 . 
\monend

\smallskip \monitem Proof of  Theorem 2. 

\noindent
From the relation  (\ref{delta-Lagrangien}), the second order adjoint equation can be written as
\moneqstar 
 \Big( {{\dd^2 \xi}\over{\dd t^2}} \Big)_j +  R^k_{\ell j m} \,  \xi_k \, \zeta^\ell \, \zeta^m 
   - {{\dd}\over{\dd t}} \Big( {{\partial \gamma}\over{\partial \zeta}} \Big)_j 
+  \Big( {{\partial \gamma}\over{\partial q}} \Big)_j
- \xi_k \, (\partial_j \psi )^k  = 0 .
\monendstar
With $\, \psi = - \partial_\ell V \, e^\ell $, we have the following calculus:

\smallskip \noindent 
$   \partial_j \psi = - \partial_j \partial_\ell V \, e^\ell  + \Gamma^s_{j \ell} \,  \partial_s V \, e^\ell
= - (\nabla^2 V)_{j \ell}  \, e^\ell \, $ 
and

\smallskip \noindent 
$ \, \xi_k \, (\partial_j \psi)^k = - M^{k \ell} \,  (\nabla^2 V)_{j \ell}  \, \xi_k = 
-  (\nabla^2 V)_{j \ell} \, \xi^\ell = -(\nabla^2 V . \,\xi)_j $.
Additionally, we establish the contraction
$ \, R^k_{\ell j m} \,  \xi_k \, \zeta^\ell \, \zeta^m = \big( R_\zeta \,.\, \xi \big)_j $. 
Then the evolution equation can be written
\moneqstar 
\Big( {{\dd^2 \xi}\over{\dd t^2}} - R_\zeta \,.\, \xi + {{\partial \gamma}\over{\partial q}}
-  {{\dd}\over{\dd t}} \big( {{\partial \gamma}\over{\partial \zeta}} \big) + \nabla^2 V . \,\xi \Big)_j = 0 . 
\monendstar 
and the relation (\ref{dynamique-second-adjoint}) is established. \hfill $\square$

\bigskip   \bigskip   \noindent {\bf \large   Conclusion}    

\smallskip \noindent
We first compared  the result for the quadratic cost function (\ref{cout-quadratique}) 
developed in paragraph~4 and the present result studied in the previous section.  
The cost function is now more general. It was written
\moneqstar
J(u) =  {1\over2} \, \int_0^T M_{k \ell} (q) \, u^k \, u^\ell \, {\rm d}t 
\monendstar
 in \cite{DFRV15}
and we write it (\ref{cost-function})
\moneqstar
J(u) = \int_0^T \gamma( q,\, \zeta ,\, u) \, \dd t 
\monendstar
in this contribution.
Nevertheless, the equations of the dynamical system take the same form:
\moneqstar
{{\dd q}\over{\dd t}} = \zeta \,,\,\,   {{\dd \zeta}\over{\dd t}} - \psi(q) =  u 
\monendstar
with the usual condition that the internal forces derive from a potential. 
With the particular cost function considered in \cite{RDR22}, 
the optimality condition take the form 
\moneqstar
\xi = u .
\monendstar
The Lagrange multiplier associated to the dynamics equation is interpreted as
a force. Then  the adjoint equation derived in 
Vall\'ee {\it et al.} \cite{Juan-these-2013,VRFG13} is exactly  a 
covariant evolution equation~(\ref{evolution-xi-old}) for the optimal force.
With the general  cost function considered in this contribution,
the  optimality condition can be written 
\moneqstar
{{\partial \gamma}\over{\partial u}} = \xi .
\monendstar

\smallskip \noindent 
The dynamics of the adjoint  variable $ \, \xi \, $ differs {\it a priori} from the one of
forces and torques~$ \, u $.  
We have explicited this condition in (\ref{dynamique-second-adjoint}).
We observe that in comparison with  (\ref{evolution-xi-old}),  
two new terms are present: 
$ \, {{\partial \gamma}\over{\partial q}} \, $ and 
$ \, -  {{\dd}\over{\dd t}} \big( {{\partial \gamma}\over{\partial \zeta}} \big) $.

\smallskip  \noindent
In this contribution, we have generalized  the cost function used for 
the Pontryagin calculus
in Riemannian geometry synthesized in \cite{DFRV15}.
The cost function is still chosen in coherence with the Riemannian geometry underlying
the natural evolution of the mechanical system.
The applications of this approach in robotics are into development
and first results are proposed in \cite{RDR22}.
The next step is the enrichement of the model with  appropriate
dissipation as fluid rubbing or dry Coulomb friction.

\bigskip  \bigskip     \noindent {\bf   \large Acknowledgments}

\noindent
The authors thank the team of International Conference Zaragoza-Pau on Applied Mathematics and Statistics
for their open-mindedness.
They thank in particular the referee who suggested a clarification
  in the presentation of this contribution.
They thank also G\'ery de Saxc\'e and Fr\'ed\'eric Boyer 
for sharing  helpful remarks after an online  presentation of this work in september 2022. 

\bigskip \bigskip   
\newpage  \noindent {\bf \large   References}    


\end{document}